# Exponential Integrators for Resistive Magnetohydrodynamics: Matrix-free Leja Interpolation and Efficient Adaptive Time Stepping

Pranab J. Deka [1] and Lukas Einkemmer [1]

[1]*Department of Mathematics*
*University of Innsbruck*
*A-6020 Innsbruck, Austria*

## ABSTRACT

We propose a novel algorithm for the temporal integration of the resistive magnetohydrodynamics (MHD) equations. The approach is based on exponential Rosenbrock schemes in combination with Leja interpolation. It naturally preserves Gauss's law for magnetism and is unencumbered by the stability constraints observed for explicit methods. Remarkable progress has been achieved in designing exponential integrators and computing the required matrix functions efficiently. However, employing them in MHD simulations of realistic physical scenarios requires a matrix-free implementation. We show how an efficient algorithm based on Leja interpolation that only uses the right-hand side of the differential equation (i.e. matrix free), can be constructed. We further demonstrate that it outperforms Krylov-based exponential integrators as well as explicit and implicit methods using test models of magnetic reconnection and the Kelvin–Helmholtz instability. Furthermore, an adaptive step-size strategy that gives excellent and predictable performance, particularly in the lenient- to intermediate-tolerance regime that is often of importance in practical applications, is employed.

## 1. INTRODUCTION

Computational methods are remarkable tools to study problems in physics, astrophysics, space physics, engineering, and several other scientific fields. They enable a better understanding of the various physical phenomena involved, and their results can be compared to experimental or observational data. This is particularly true for plasma physics and astrophysics. They are vital where the experimental or observational measurements are prone to large errors. The macroscopic behaviour of a plasma or any conducting fluid can be efficiently described by coupling Maxwell's equations with hydrodynamics (the Euler equations or the Navier-Stokes equations), thereby giving rise to equations of MHD. MHD is based on the principle that a dynamical magnetic field can induce current in a conducting fluid, which can, in turn, change the magnetic field of the system. There are a plethora of codes that solve the MHD equations. A non-exhaustive list comprises of ATHENA (Stone et al. 2008), CASTRO (Almgren et al. 2010), CRONOS (Kissmann et al. 2018), JOREK (Huysmans & Czarny 2007), NIRVANA (Ziegler 2008), RACOON (Dreher & Grauer 2005), RAMSES (Teyssier 2002; Fromang et al. 2006), PLUTO (Mignone et al. 2007), and many more. The focus of most of these MHD codes has been on adaptive meshes and relativistic treatment of particles at very high energies. They are usually tailored to model specific physical or astrophysical scenarios. Running such simulations can often take days or even weeks and, in some cases, require the usage of large supercomputers. It is worth noting that most of the aforementioned codes use a second (or in some cases, third) order explicit (e.g. RK2, RK3) or implicit integrators (e.g. Crank-Nicolson). The use of explicit integrators restricts the step size to the Courant-Friedrich-Lewy (CFL) limit, for both constant and adaptive step sizes. Implicit solvers usually do not suffer from such restrictions but are often computationally demanding and difficult to implement.

In this paper, we propose a novel algorithm for the temporal integration of the MHD equations. This algorithm is based on exponential Rosenbrock methods that use Leja interpolation to compute the corresponding matrix functions. Our implementation is completely matrix-free, i.e. similar to explicit methods where only the evaluations of the right-hand side are required and not restricted by the CFL condition that explicit schemes have to adhere to. The divergence-free constraint on the magnetic field is preserved inherently (i.e. Gauss's law for magnetism is satisfied without the need to perform divergence cleaning). We investigate the performance of the proposed algorithm and show that it outperforms explicit and exponential integrators (based on Krylov iterations) in a variety of circumstances.

Exponential Rosenbrock (EXPRB) methods are a class of numerical time integrators that linearize the underlying differential equation and treat the linear part 'exactly'. This is why they preserve many physical features of the system, e.g. all linear invariants of the original system and equilibria are automatically preserved in exponential integrators (Theorem 2 in Einkemmer & Ostermann (2015)). Similar to implicit and semi-implicit methods, there is a need to efficiently compute certain matrix functions for the exponential Rosenbrock methods (more details in Sec. 2). The Krylov subspace projection algorithm (Arnoldi 1951; Van Der Vorst 1987) proves to be an excellent tool for this purpose if the matrices obtained from the discretization of the partial differential equations (PDEs) are large and sparse, which is often the case in most problems in physics. Einkemmer et al. (2017) studied the performance of the exponential integrators using the Exponential Propagation Integrators Collection (EPIC, Tokman (2014)) package, which implements the aforementioned Krylov iteration, applied to the MHD equations. It was found that exponential Rosenbrock-type schemes show improved performance compared to traditionally used implicit schemes (specifically, the commonly used CVODE library) for a range of configurations, i.e. the user-defined tolerance and the error incurred for adaptive step sizes, different values of constant step sizes, and the parameters of viscous MHD ($\mu$, $\eta$, $\kappa$, see Sec. 3) for two different MHD problems.

Using polynomial interpolation at Leja points to compute the matrix functions offers an attractive alternative to the Krylov subspace methods. Such an approach only requires matrix-vector products and is thus computationally cheaper. This is particularly true for massively parallel architectures (such as graphical processing units, GPUs) that are the future of high-performance computing. These favourable properties of Leja interpolation for the implementation of exponential Rosenbrock methods is well established in the literature, see (Bergamaschi et al. 2006; Caliari et al. 2014; Auer et al. 2018). However, in most of these works, relatively simple test examples are considered in which the matrix of the linearization is given explicitly by an analytical formula. In the case of the (highly nonlinear) MHD equations, assembling this matrix would be prohibitively costly.

Employing Leja interpolation requires an estimate of the largest eigenvalue. It is straightforward to obtain such an estimate when the matrix is explicitly assembled. However, in the case of MHD problems, it is of utmost importance to use a matrix-free formulation, i.e. only evaluations of the right-hand side are needed in the time-integration scheme. This is essential for computational efficacy and to readily integrate such a method into existing codes. In this paper, we demonstrate that the information on the spectrum of the linearization can be obtained with negligible additional cost by performing power iterations. This, to the best of our knowledge, is the first matrix-free Leja-based exponential integrator described in the literature.

Another crucial aspect of making a numerical method attractive for MHD is to free the user from selecting an appropriate time step size. For explicit methods, this is often done by simply operating well below the CFL limit. Exponential integrators, however, due to the exact treatment of the linearization, can take large time steps (often significantly larger than even fully implicit methods). However, selecting a very large time step is not always beneficial from the computational point of view. The reason for this is that taking larger time steps increases stiffness, which can result in the Leja interpolation (or Krylov iteration) to converging slower. It is an interesting observation that selecting the largest possible time step based on the user-specified tolerance is not necessarily optimal from the performance standpoint. Most traditional step-size controllers operate under the assumption that the computational cost is independent of the step size. This is valid only for explicit methods and not for iterative solvers, which is the case for implicit methods as well as exponential integrators. However, for the latter, the issue is more severe as, in principle, even larger time steps are possible. An adaptive step-size controller that takes this into account has been proposed for implicit methods (Einkemmer 2018). More recently, the authors have extended it to exponential integrators (Deka & Einkemmer 2021). We use this approach in the present work. A further benefit of this step-size controller is that the error obtained in the solution more directly corresponds to the specified tolerance (see Sec. 5).

We study the performance of the proposed approach in the context of magnetic reconnection and the Kelvin–Helmholtz instability (KHI). We demonstrate that our algorithm is particularly efficient in the lenient- to medium-tolerance range (which is of interest in many applications), where significant computational gains can be observed in virtually all configurations (up to a factor of three when comparing Leja interpolation to Krylov iteration). We also conduct comparisons to explicit integrators and obtain still larger gains in performance.

This manuscript is structured as follows: In Sec. 2, we give a gist of the exponential integrators, present the relevant schemes for this study, and provide a mathematical analysis of the step-size controller under consideration. The MHD equations are described in Sec. 3. The performance of the step-size controllers with a fourth-order scheme is analyzed for the KHI in Sec. 4. Here, we also contrast the performance of exponential integrators with explicit methods. The performance of a fifth-order integrator is examined for magnetic reconnection, and a comparison of the Leja and Krylov methods is articulated in Sec. 5. Finally, we conclude in Sec. 6.

## 2. EXPONENTIAL INTEGRATORS: EFFICIENT MATRIX-FREE IMPLEMENTATION

Systems of differential equations where the characteristic time scales of the relevant physical processes vary drastically can be classified as *stiff* systems. Explicit time integrators are forced to choose very small step sizes to preserve the numerical stability. Furthermore, diffusion dominated problems are severely constricted by the CFL limit. Although implicit integrators are better suited to preserve the stability of the solution, one has to solve a large system of nonlinear equations, by means of Newton iterations, at every time step. Typical methods of solving the resulting system of large linear systems include conjugate gradient, GMRES, biconjugate gradient, etc. However, one has to incur large computational expenses for solving such large systems of linear equations.

Exponential integrators are a class of temporal integrators suitable for stiff systems described by one or a system of differential equations. They have excellent stability properties, and they require the solution of the matrix exponential (or related functions) instead of solutions to large linear systems. Consequently, it becomes necessary to solve the matrix exponential in an efficient and cost-effective way (more in Sec. 2.3). Here, we present the basic idea of exponential integrators; for further details on the technical aspects, we refer the reader to the the review article by Hochbruck & Ostermann (2010). Let us consider the following initial value problem:

$$\frac{\partial u}{\partial t} = f(u, t), \quad u(t = 0) = u^0, \quad (1)$$

where $f(u, t)$ is some nonlinear function of $u$. In the case of MHD problems, this equation results from the discretization of the corresponding PDE in space (i.e. the method of lines). The general idea of exponential integrators is to write Eq. 1 as follows

$$\frac{\partial u}{\partial t} = \mathcal{A} u + g(u), \quad (2)$$

where $\mathcal{A}$ is a matrix and $g(u)$ is a non-linear remainder. The simplest exponential integrator is the exponential Euler method given by

$$u^{n+1} = u^n + \varphi_1(\mathcal{A}\Delta t) f(u^n) \Delta t.$$

This is, in general, a first-order scheme that requires one evaluation of the right-hand side of Eq. 1 and the action of the matrix function $\varphi_1(\mathcal{A}\Delta t)$ to the vector $f(u^n)\Delta t$. The $\varphi(z)$ functions are defined recursively as follows:

$$\varphi_{l+1}(z) = \frac{1}{z}\left(\varphi_l(z) - \frac{1}{l!}\right), \quad l \geq 1$$

with

$$\varphi_0(z) = e^z$$

corresponding to the matrix exponential. In the case of linear PDEs (that is, if $f(u, t) = \mathcal{A} u$, and $g(u) = 0$), one can obtain the exact solution of the system, given by

$$u^{n+1} = u^n \exp(\mathcal{A}\Delta t).$$

In some problems, for an appropriate choice of $\mathcal{A}$, the action of the matrix function can be computed directly (using FFT techniques). This approach has been used, among others, for kinetic equations (Crouseilles et al. 2018; Caliari et al. 2021), water waves (Klein & Roidot 2011; Caliari et al. 2021), sonic boom propagation (Einkemmer et al. 2021), and the simulation of Bose–Einstein condensates (Besse et al. 2017). However, due to the fully nonlinear nature and the fundamentally complex structure of the MHD equations, such an approach is not feasible. Therefore, an iterative scheme has to be used. In the following, we will discuss the fourth- and fifth-order embedded exponential integrators that we will use in this work. Details on the order conditions and the construction of such methods can be found in Hochbruck & Ostermann (2010).

### 2.1. *Exponential Rosenbrock Schemes*

EXPRB schemes are a sub-class of exponential integrators that use the Jacobian of $f(u)$ to define the linear part $\mathcal{A}$. That is, at each time step, we linearize at $u^n$ to obtain

$$\frac{\partial u}{\partial t} = \mathcal{J}(u^n) u + \mathcal{F}(u),$$

where $\mathcal{J}(u)$ is the Jacobian of $f(u)$ and $\mathcal{F}(u) = f(u) - \mathcal{J}(u)$ is the non-linear remainder. The *Rosenbrock-Euler* method, the simplest of the exponential Rosenbrock methods, is given by

$$u^{n+1} = u^n + \varphi_1(\mathcal{J}(u^n)\Delta t) f(u^n) \Delta t. \quad (3)$$

We duly note that due to the use of the Jacobian, this is, in fact, a second-order method.

In the present work, we consider the fourth-order EXPRB43 scheme embedded with a third-order error estimate. The equations of the EXPRB43 scheme (Hochbruck et al. 2009; Hochbruck & Ostermann 2010) are given by

$$a^n = u^n + \frac{1}{2}\varphi_1\left(\frac{1}{2}\mathcal{J}(u^n)\Delta t\right) f(u^n)\Delta t$$

$$\begin{aligned} b^n = u^n &+ \varphi_1\left(\mathcal{J}(u^n)\Delta t\right) f(u^n)\Delta t \\ &+ \varphi_1\left(\mathcal{J}(u^n)\Delta t\right)(\mathcal{F}(a^n) - \mathcal{F}(u^n))\Delta t \end{aligned}$$

$$\begin{aligned} \hat{u}^{n+1} = u^n &+ \varphi_1\left(\mathcal{J}(u^n)\Delta t\right) f(u^n)\Delta t \\ &+ \varphi_3(\mathcal{J}(u^n)\Delta t)(-14\mathcal{F}(u^n) + 16\mathcal{F}(a^n) - 2\mathcal{F}(b^n))\Delta t \end{aligned} \quad (4)$$

$$u^{n+1} = u^n + \varphi_1\left(\mathcal{J}(u^n)\Delta t\right) f(u^n)\Delta t$$
$$+ \varphi_3(\mathcal{J}(u^n)\Delta t)(-14\mathcal{F}(u^n) + 16\mathcal{F}(a^n) - 2\mathcal{F}(b^n))\Delta t$$
$$+ \varphi_4(\mathcal{J}(u^n)\Delta t)(36\mathcal{F}(u^n) - 48\mathcal{F}(a^n) + 12\mathcal{F}(b^n))\Delta t, \quad (5)$$

where $a_n$ and $b_n$ are the internal stages, and Eq. 4 and 5 correspond to the third- and the fourth-order solutions, respectively. The difference between the third- and fourth-order solutions generates a third-order error estimate. In Secs. 4 and 5, we will see that the fourth-order EXPRB43 method presents an excellent trade-off between accuracy and the computational cost.

In addition to the fourth-order scheme, we test the performance of a 4-stage fifth-order method, EXPRB54s4 (Luan & Ostermann 2014), embedded with a fourth-order error estimator. This is to facilitate a reasonable comparison of the EXPRB methods with the fifth-order EPIRK5P1 scheme (see following subsection). This scheme reads as follows:

$$a^n = u^n + \frac{1}{4}\varphi_1\left(\frac{1}{4}\mathcal{J}(u^n)\Delta t\right) f(u^n)\Delta t$$
$$b^n = u^n + \frac{1}{2}\varphi_1\left(\frac{1}{2}\mathcal{J}(u^n)\Delta t\right) f(u^n)\Delta t$$
$$+ 4\varphi_3\left(\frac{1}{2}\mathcal{J}(u^n)\Delta t\right)(\mathcal{F}(a^n) - \mathcal{F}(u^n))\Delta t$$
$$c^n = u^n + \frac{9}{10}\varphi_1\left(\frac{9}{10}\mathcal{J}(u^n)\Delta t\right) f(u^n)\Delta t$$
$$+ \frac{729}{125}\varphi_3\left(\frac{9}{10}\mathcal{J}(u^n)\Delta t\right)(\mathcal{F}(b^n) - \mathcal{F}(u^n))\Delta t$$

$$\hat{u}^{n+1} = u^n + \varphi_1\left(\mathcal{J}(u^n)\Delta t\right) f(u^n)\Delta t$$
$$+ \varphi_3(\mathcal{J}(u^n)\Delta t)\mathcal{F}_1\Delta t + \varphi_4(\mathcal{J}(u^n)\Delta t)\mathcal{F}_2\Delta t \quad (6)$$

$$u^{n+1} = u^n + \varphi_1\left(\mathcal{J}(u^n)\Delta t\right) f(u^n)\Delta t$$
$$+ \varphi_3(\mathcal{J}(u^n)\Delta t)\mathcal{F}_3\Delta t + \varphi_4(\mathcal{J}(u^n)\Delta t)\mathcal{F}_4\Delta t \quad (7)$$

where $\mathcal{F}_1 = -56\mathcal{F}(u^n) + 64\mathcal{F}(a^n) - 8\mathcal{F}(b^n)$, $\mathcal{F}_2 = 80\mathcal{F}(u^n) - 60\mathcal{F}(a^n) - \frac{285}{8}\mathcal{F}(b^n) + \frac{125}{8}\mathcal{F}(c^n)$, $\mathcal{F}_3 = -\frac{1208}{81}\mathcal{F}(u^n) + 18\mathcal{F}(b^n) - \frac{250}{81}\mathcal{F}(c^n)$, and $\mathcal{F}_4 = \frac{1120}{27}\mathcal{F}(u^n) - 60\mathcal{F}(b^n) + \frac{500}{27}\mathcal{F}(c^n)$. Eqs. 6 and 7 give the fourth- and fifth-order solutions, respectively, and the difference between them gives the fourth-order error estimate.

### 2.2. Exponential Propagation Iterative Runge–Kutta Schemes

Exponential Propagation Iterative Runge–Kutta (EPIRK) schemes (Tokman 2006) constitute another class of exponential integrators. These schemes are engineered specifically to enable reusing and recycling of the constructed Krylov space for multiple matrix functions, i.e. they are "Krylov-friendly schemes." In this work, we consider the fifth-order EPIRK5P1 scheme (Tokman et al. 2012) embedded with a fourth-order error estimator:

$$a^n = u^n + a_{11}\,\varphi_1\left(g_{11}\,\mathcal{J}(u^n)\Delta t\right) f(u^n)\Delta t$$
$$b^n = u^n + a_{21}\,\varphi_1\left(g_{21}\,\mathcal{J}(u^n)\Delta t\right) f(u^n)\Delta t$$
$$+ a_{22}\,\varphi_1\left(g_{22}\,\mathcal{J}(u^n)\Delta t\right)(\mathcal{F}(a^n) - \mathcal{F}(u^n))\Delta t$$

$$u^{n+1} = u^n + b_1\,\varphi_1\left(g_{31}\,\mathcal{J}(u^n)\Delta t\right) f(u^n)\Delta t$$
$$+ b_2\,\varphi_1\left(g_{32}\,\mathcal{J}(u^n)\Delta t\right)(\mathcal{F}(a^n) - \mathcal{F}(u^n))\Delta t$$
$$+ b_3\,\varphi_3\left(g_{33}\,\mathcal{J}(u^n)\Delta t\right)(\mathcal{F}(u^n) - 2\mathcal{F}(a^n) + \mathcal{F}(b^n))\Delta t \quad (8)$$

The coefficients corresponding to this scheme can be found in Table 1. Here, $a_n$ and $b_n$ correspond to the internal stages, and Eq. 8 gives the fifth-order solution. The embedded fourth-order solution can be obtained by simply setting $g_{32} = 0.5$ and $g_{33} = 1.0$ in Eq. 8 (Loffeld & Tokman 2013). The difference between these two solutions gives an error estimate of order four.

### 2.3. Polynomial interpolation at Leja points

The most computationally challenging part of an exponential integrator is the action of the matrix function on a vector. Although traditional methods like the Padé approximation or Taylor expansion work well for small matrices, they are computationally prohibitive for large matrices. Krylov subspace methods and polynomial interpolation can be viable alternatives for large-scale problems (Bergamaschi et al. 2006; Caliari et al. 2014). In this work, we primarily focus on evaluating the matrix functions by interpolating them as polynomials on a suitable subset consisting of its spectrum. Although Chebyshev points are widely used for this purpose, they suffer from the drawback that interpolation of a function at 'm + 1' Chebyshev points requires the reevaluation of the function at all points (no reuse of previously computed values is possible). This can lead to large computational cost if a large number of nodes are needed to efficiently interpolate the given function. Leja points, however, can be defined recursively as

$$z_m \in \arg\max \prod_{j=0}^{m-1} |z - z_j|$$

for given $z_0$, $z_m$ is the $m^{\text{th}}$ Leja point, and $z \in K \subseteq \mathbb{C}$, $\mathbb{C}$ being the set of complex numbers. Here, the evaluation of a function at 'm + 1' points, if already interpolated at 'm' Leja points, needs only one extra function evaluation. This makes the use of Leja points attractive from a computational point of view. For further details

$$\begin{bmatrix} a_{11} & & \\ a_{21} & a_{22} & \\ b_1 & b_2 & b_3 \end{bmatrix} = \begin{bmatrix} 0.35129592695058193092 & & \\ 0.84405472011657126298 & 1.6905891609568963624 & \\ 1.0 & 1.2727127317356892397 & 2.2714599265422622275 \end{bmatrix}$$

$$\begin{bmatrix} g_{11} & & \\ g_{21} & g_{22} & \\ g_{31} & g_{32} & g_{33} \end{bmatrix} = \begin{bmatrix} 0.35129592695058193092 & & \\ 0.84405472011657126298 & 0.5 & \\ 1.0 & 0.71111095364366870359 & 0.62378111953371494809 \end{bmatrix}$$

**Table 1.** Coefficients of EPIRK5P1

on the polynomial interpolation method at Leja points, we refer the reader to the manuscripts by Caliari et al. (2004, 2007, 2014).

To compute the Jacobian ($\mathcal{J}$) in our matrix-free implementation, we evaluate the right-hand side using forward finite differences. The implementation of Leja interpolation requires an estimate of the largest eigenvalue of the Jacobian (henceforth denoted by $\alpha$). Whilst such an estimate is easy to compute if the matrix is assembled, here, we obtain an estimate by using power iterations (subject to a safety factor). At first glance, this procedure might seem expensive, but this problem can be dealt with simply and efficiently (Sec. 5.2). This is the primary difference of the present algorithm with that of Deka & Einkemmer (2021). As in Caliari et al. (2014), we shift the matrix around the midpoint $q = \frac{1}{2}\alpha$ and scale it by $\theta = \frac{1}{4}\alpha$ (in our case, the Leja points are defined on the interval $[-2, 2]$).

One can now interpolate the function $\varphi(q+\theta\xi)$ for $\xi \in [-2, 2]$. The $m^{\text{th}}$ term of the interpolation polynomial $p(z)$ is given by

$$p_m(z) = p_{m-1}(z) + d_m\, y_{m-1}(z),$$
$$y_m(z) = y_{m-1}(z) \times \left(\frac{z-q}{\theta} - \xi_m\right),$$

where $d_i$ corresponds to the divided differences of the relevant $\varphi$ functions. This can easily be turned into a scheme for computing the action of matrix polynomials simply by replacing $z$ with the appropriate action of the Jacobian.

We like to point out one caveat in using the Leja interpolation method. The convergence of a multi-scheme exponential integrator depends on the convergence of each of the internal stages to the user-specified tolerance. This convergence depends on the norm of the function being interpolated as a polynomial and the magnitude of $\Delta t$. If both of these parameters happen to be too large, the polynomial may not converge; one needs to reject the step size and start with a smaller one as dictated by the traditional controller (described in the following section).

### 2.4. Adaptive Step-Size Controller

Adaptive step-size controllers selectively modify the step size to match the accuracy requirements set by the user. The widely used traditional step size controller uses an error estimate from the previous time step. The new step size is chosen to be the largest possible step size subject to the user-specified tolerance (with a safety factor). This approach implicitly assumes that the computational cost is independent of the step size. This is reasonable for explicit integrators or implicit integrator that uses direct solvers. The traditional step-size controller is given by

$$\Delta t^{n+1} = \Delta t^n \times \left(\frac{\text{tol}}{e^n}\right)^{1/(p+1)},$$

where $\Delta t^n$ and $\Delta t^{n+1}$ correspond to the step sizes at the present ($n^{\text{th}}$) and the next time step (($n+1$)$^{\text{th}}$), respectively. The error incurred at the $n^{\text{th}}$ time step is denoted by $e^n$, $p$ is the order of the integration scheme, and tol is the user-specified tolerance. The mathematical analysis of these kinds of step-size controllers can be found in Gustafsson et al. (1988), Gustafsson (1994), Söderlind (2002), and Söderlind (2006). Variations of this step size controller are used in several time-integration software packages, e.g. `RADAU5` (Hairer & Wanner 1996) and `CVODE` (Eckert et al. 2004; Hindmarsh & Serban 2016) to name a few.

The assumption that the cost of a time step is independent of the step size is not valid for iterative schemes. This is true for both implicit and exponential integrators. In this case, the computational cost incurred depends strongly and in a nonlinear fashion on the step size. Thus, using the largest possible step size as dictated by the traditional controller might not be advantageous. This is the principle behind the step-size controller proposed in Einkemmer (2018), which has been developed to minimize the computational cost. Consequently, the step sizes chosen by this controller might be significantly smaller than the ones yielded by the traditional controller (which chooses the largest possible step size based on the given tolerance).

The computational cost ($c$) per unit time step size is given by

$$c^n = \frac{i^n}{\Delta t^n},$$

where $i^n$ is the computational runtime (or some proxy of the computational runtime). The goal is to minimize the cost: $c^n \longrightarrow \min$. We work with the logarithm of the step size $T^n = \ln \Delta t^n$ and the cost $C^n(T^n) = \ln c^n(\Delta t^n)$. Using gradient descent one obtains

$$T^{n+1} = T^n - \Gamma \nabla C^n(T^n),$$

where $\Gamma$ is the learning rate. The gradient is then approximated by finite differences

$$\nabla C^n(T^n) \approx \frac{C^n(T^n) - C^n(T^{n-1})}{T^n - T^{n-1}}.$$

It is important to appreciate the difference between $C^n(T^n)$, the cost of a time step with size $\Delta t^n$ at the $n^{\text{th}}$ time step, and $C^n(T^{n-1})$, the cost of the same step size $\Delta t^n$ at the $(n-1)^{\text{th}}$ time step. The temporal integration of a problem yields $C^n(T^n)$ and *not* $C^n(T^{n-1})$. Taking this into account, we can approximate the gradient as

$$\begin{aligned}\nabla C^n(T^n) &\approx \frac{C^n(T^n) - C^n(T^{n-1})}{T^n - T^{n-1}} \\ &= \frac{C^n(T^n) - C^{n-1}(T^{n-1})}{T^n - T^{n-1}} \\ &\quad + \frac{C^{n-1}(T^{n-1}) - C^n(T^{n-1})}{T^n - T^{n-1}} \\ &\approx \frac{C^n(T^n) - C^{n-1}(T^{n-1})}{T^n - T^{n-1}},\end{aligned}$$

assuming that $C^n$ varies slowly as a function of $n$. This gives

$$T^{n+1} = T^n - \Gamma \frac{C^n(T^n) - C^{n-1}(T^{n-1})}{T^n - T^{n-1}}.$$

Taking exponential on both sides of the equation, we get

$$\Delta t^{n+1} = \Delta t^n \exp(-\Gamma \Delta), \qquad \Delta = \frac{\ln c^n - \ln c^{n-1}}{\ln \Delta t^n - \ln \Delta t^{n-1}}.$$

The free parameter $\Gamma$ can be a function of either $c$, or $\Delta t$, or both. Choosing $\Gamma$ to be a constant, despite being the simplest choice, has the disadvantage that it may either lead to minute (ineffective) changes or prohibitively large changes in $\Delta t$. This is why the new step size is computed as follows:

$$\Delta t^{n+1} = \Delta t^n \times \begin{cases} \lambda & \text{if } 1 \leq s < \lambda, \\ \delta & \text{if } \delta \leq s < 1, \\ s := \exp(-\alpha \tanh(\beta \Delta)) & \text{otherwise.} \end{cases}$$

$\alpha$ is a constraint on the largest possible change in stepsize, given by $\exp(\pm\alpha)$. $\beta$ is the response of the step size controller to the change in cost. $\lambda$ and $\delta$ have been incorporated to ensure that the step size changes by at least $\lambda \cdot \Delta t$ or $\delta \cdot \Delta t$ depending on whether $\Delta t$ needs to be increased or decreased to minimize the cost.

The parameters of the step-size controller have been optimized, in Einkemmer (2018), to incur the least possible computational cost, and we refer the reader to this work for more details. The numerical values of the parameters of the proposed controller are as follows: $\alpha = 0.65241444$, $\beta = 0.26862269$, $\lambda = 1.37412002$, and $\delta = 0.64446017$.

It is worth noting that the proposed step-size controller solely takes into account the effect of minimizing the computational cost. This necessitates that the error incurred at each time step is considered by some other means, i.e. the traditional controller. Therefore, we choose the minimum of the two step sizes yielded by the two controllers:

$$\Delta t^{n+1} = \min\left(\Delta t^{n+1}_{\text{proposed}}, \Delta t^{n+1}_{\text{traditional}}\right).$$

The initial study (Einkemmer 2018) of this proposed step-size controller was performed using implicit integrators (namely, Crank–Nicolson and implicit Runge–Kutta methods) for a range of linear and nonlinear PDEs using the iterative Krylov subspace method. In Deka & Einkemmer (2021), we have shown that the proposed controller works equally well with exponential integrators employing the method of polynomial interpolation at Leja points. In both these studies, one can see that the proposed step-size controller has superior performance compared to that of the traditional controller for a vast majority of the considered configurations, thereby validating the hypothesis that multiple small step sizes do indeed incur less computational cost than a single large step size.

## 3. THE MHD EQUATIONS

We introduce the equations of fluid dynamics in brief. For details on the kinetic and the fluid theory, derivation of the equations, its applicability, and the limits of validation, we refer the reader to Chen (1984); Goedbloed et al. (2010); Chiuderi & Velli (2010). The kinetic theory of plasma considers plasma to be an ensemble of charged particles, namely electrons and ions. The particle distribution function is evolved under the effect of the self-consistent electromagnetic fields. One needs to solve the Boltzmann equation (or the Vlasov equation, if particle collisions are neglected) to obtain the kinetic description. The fluid model of a plasma, however, is based on its description of the macroscopic quantities (i.e. density, velocity, and energy). The fluid approximation holds true if the length scales under consideration are much larger than the plasma skin length and the gyroradius. The equations for the fluid description can be obtained by taking velocity moments of the Boltzmann equation. In situations where the fluid approximation is valid, the equations of MHD (Maxwell's equations coupled to the fluid equations) are most commonly used to describe plasma dynamics. The MHD approximation further assumes that the ions and electrons have a common velocity.

In a real fluid, the magnetic field lines can move through the fluid, which inherently introduces resistivity and viscosity. Furthermore, if the collisional timescales are comparable to the timescales under consideration, one needs to consider dissipative effects (thermal conductivity) due to collisions. The combined effects of resistivity, viscosity, and thermal conductivity give rise to resistive MHD, the equations of which, in conservative form, read

$$\frac{\partial U}{\partial t} + \nabla \cdot F_{\text{ideal}}(U) = \nabla \cdot F_{\text{diff}}(U) \qquad (9)$$

with

$$F_{\text{ideal}}(U) = \begin{bmatrix} \rho \boldsymbol{v} \\ \rho(\boldsymbol{v} \otimes \boldsymbol{v}) + \left(P + \frac{B^2}{2\mu_0}\right) I - \frac{1}{\mu_0}(\boldsymbol{B} \otimes \boldsymbol{B}) \\ (\boldsymbol{v} \otimes \boldsymbol{B}) - (\boldsymbol{B} \otimes \boldsymbol{v}) \\ \left(E + P + \frac{B^2}{2\mu_0}\right) \boldsymbol{v} - \frac{1}{\mu_0} \boldsymbol{B}(\boldsymbol{B} \cdot \boldsymbol{v}) \end{bmatrix}$$

and

$$F_{\text{diff}}(U) = \begin{bmatrix} 0 \\ \mu \boldsymbol{\tau} \\ \eta \left(\nabla \otimes \boldsymbol{B} - (\nabla \otimes \boldsymbol{B})^T\right) \\ \mu \boldsymbol{\tau} \cdot \boldsymbol{v} + \mu \kappa \frac{\gamma}{\gamma - 1} \nabla T + \eta \left(\frac{1}{2} \nabla(\boldsymbol{B} \cdot \boldsymbol{B}) - \boldsymbol{B}(\nabla \otimes \boldsymbol{B})^T\right) \end{bmatrix}$$

$$\boldsymbol{\tau} = \nabla \boldsymbol{v} + \nabla \boldsymbol{v}^T - \frac{2}{3} \nabla \cdot \boldsymbol{v} I,$$

Here, $\rho$ corresponds to the mass density, $\boldsymbol{v}$ is the fluid velocity, $\boldsymbol{B}$ is the magnetic field, $E$ is the total energy density of the system, $I$ is the identity matrix, and $\otimes$ denotes the tensor product. This set of equations can be closed by determining the pressure ($P$) from the isentropic relation

$$\left(\frac{P}{\rho}\right)^\gamma = \text{constant} \qquad (10)$$

where we have assumed an adiabatic system. For a single-fluid approximation, we have specific heat, $\gamma = 5/3$. $F_{\text{diff}}(U)$ takes into account the dissipative effects in a plasma. Here, we have assumed constant viscosity $\mu$, resistivity $\eta$, and thermal conductivity $\kappa$. One can define the Reynold's number $Re = \rho_0 v_A l_0 / \mu_0$, Lundquist number $S = \mu_0 v_A l_0 / \eta_0$, and Prandtl number $Pr = C_p \mu_0 / \kappa_0$ for characteristic length scales $l_0$, characteristic density $\rho_0$, characteristic thermal conductivity $\kappa_0$, Alfven velocity $v_A$, permeability of free space $\mu_0$, and specific heat at constant pressure $C_P$. All numerical simulations are conducted in the dimensionless formulation. The relevant parameters (for resistive MHD) are given in terms of dimensionless viscosity $\mu = Re^{-1}$, dimensionless resistivity $\eta = S^{-1}$, and dimensionless thermal conductivity $\kappa = Pr^{-1}$.

### 3.1. *Implementation*

We use the `cppmhd` code (Einkemmer 2016), a 2.5 dimensional MHD code (velocity and magnetic field are three dimensional, and all other physical parameters are two-dimensional vectors) written in C++ with effective parallel implementation employing the Message Passing Interface library for large-scale problems. This code has been shown to produce results in Einkemmer et al. (2017) similar to those in Reynolds et al. (2006, 2010). As the Jacobian ($\mathcal{J}$) of the MHD equations has been computed using finite differences, the discrete solenoidal property, i.e.

$$\frac{U_{i+1,j}^n - U_{i-1,j}^n}{2\Delta x} + \frac{U_{i,j+1}^n - U_{i,j-1}^n}{2\Delta y} = 0, \qquad (11)$$

holds true for $\mathcal{J}$ as well as $\mathcal{F}(U) - \mathcal{J}$. Theorem 2 in Einkemmer & Ostermann (2015) states that any exponential Runge–Kutta method satisfies this property. Since the computation of $\varphi(z)\,\boldsymbol{r}$ (for some vector $\boldsymbol{r}$) as a polynomial satisfies the linear constraint, a conserved quantity (here, $\nabla \cdot \boldsymbol{B} = 0$) at the initial time is also conserved at any later time. This constraint on the divergence of the magnetic field is satisfied and validated for exponential integrators (using Krylov methods) in Einkemmer et al. (2017) and (using Leja interpolation) in Sec. 5.2.

`cppmhd` provides a function to evaluate the right-hand side of the MHD equations. This is the only problem-specific function for the matrix-free implementation of the EXPRB scheme/Leja interpolation. The proposed step-size controller has also been incorporated into some of the relevant integrators in the `EPIC` library for the purposes of this work.

## 4. KELVIN-HELMHOLTZ INSTABILITY

In this and the following section, we consider two model problems based on the resistive MHD equations - the KHI and magnetic reconnection, respectively. We compare the performance of the Krylov subspace algorithm with that of the Leja interpolation method and the traditional controller with the proposed step-size controller for exponential integrators.

The KHI develops when there arises a velocity gradient at the interface of two adjacent fluids or when a fluid has continuous velocity shear. Owing to its ubiquity in nature, the KHI is an active area of research (Borse et al. 2021; Schilling 2021; Yang et al. 2021). In this work, we have chosen a simple example (similar to the one in Reynolds et al. (2010)) to simulate the KHI.

### 4.1. *Simulation setup*

The simulation domain is chosen to be $X \times Y = [-1.25, 1.25] \times [-0.5, 0.5]$. We consider a sharp velocity gradient along the Y-axis given by Eq. 12. A perturbation ($v_{\text{pert}}$) is introduced along the X-axis on $\boldsymbol{v}$ (Eq. 13). The density has been initialized to unity. Magnetic field and pressure have been initialized with constant values. The initial conditions for KHI are summarized in Table 2. Periodic boundary conditions are considered along both X and Y directions.

**Table 2.** Initial conditions for Kelvin-Helmholtz instability

| Parameters | Value |
|---|---|
| $\epsilon_x, \epsilon_y$ | 0.1 |
| $\omega_x, \omega_y$ | 2 |
| $\xi$ | 0.1 |
| $\rho$ (density) | 1.0 |
| $P$ (pressure) | 0.25 |
| $(B_x, B_y, B_z)$ | (0.1, 0, 10) |

$$\boldsymbol{v} = \begin{bmatrix} v_0 \, \tanh(y/\xi) + v_{\text{pert}} \\ 0 \\ 0 \end{bmatrix} \quad (12)$$

$$v_{\text{pert}} = \epsilon_x \cos\left(\frac{2\pi\omega_x x}{L_x}\right) + \epsilon_y \sin\left(\frac{\pi(2\omega_y - 1)y}{L_y}\right) \quad (13)$$

We use the EXPRB43 scheme for this example, and we consider the following configurations:
**Case I:** $N_x \times N_y = 512 \times 512$; $T_f = 1.0$
**Case II:** $N_x \times N_y = 800 \times 800$; $T_f = 0.3$
**Case III:** $N_x \times N_y = 128 \times 128$; $T_f = 2.0$
**Case IV:** $N_x \times N_y = 256 \times 256$; $T_f = 1.0$
where $N_x, N_y$ correspond to the number of grid points along the X and Y directions, respectively, and $T_f$ is the final simulation time. The values of $\mu$, $\eta$, and $\kappa$ are chosen to be 0.25, $10^{-2}$, and $10^{-4}$ for cases I and II (i.e. diffusion is large) and $10^{-4}$, $10^{-4}$, and $10^{-4}$ for cases III and IV (i.e. diffusion is small), respectively.

### 4.2. Results and Discussion

Fig. 1 shows the evolution of velocity with time. Figs. 2, 3, 4, and 5 compare the performance of the Leja interpolation method (dashed-dotted lines) with that of the Krylov subspace algorithm (dotted lines) for cases I, II, III, and IV, respectively. The computational expense, measured in terms of time elapsed during the computations, is plotted as a function of the user-defined tolerance (left) and the l2 norm of global error (right). The reference solution (to compute the global error) is determined with a user-specified tolerance of $10^{-11}$.

We note that there are three important quantities for a reasonable comparison of the results – the user-defined tolerance, the error incurred, and the computational runtime. The tolerance, being an input parameter, is fully under the control of the user, whilst the error incurred and the runtime are output parameters. The error incurred during a simulation depends on the fundamental scheme, the integrating method, the implementation, and the problem under consideration and its specifications, including the initial and boundary conditions. A comparison of the different methods based on the error incurred inherently takes into account all these factors, which is why it is unlikely that for a given value of tolerance, different methods or integrators will yield the same error. Since prior to conducting a simulation, it is impossible to determine the error incurred, and it is only the tolerance that can be specified by the user:comparing the computational runtime as a function of the user-defined tolerance (input) is of much practical interest (assuming, of course, that the scheme succeeds in achieving an error that lies below the specified tolerance). In this work, we propose the use of Leja-based exponential integrators primarily for the purpose of speeding up (large-scale) simulations for a given tolerance.

First, we compare the performance of the two step-size controllers. One can see that the proposed step-size controller shows an improvement in performance over the traditional controller for all values of the user-defined tolerance for both the Leja and Krylov methods (cases I and II). Moreover, it is to be noted that whilst using the traditional controller (i.e. large step sizes) for the Krylov subspace method, the simulations converge only for tol = $10^{-5}$ and $10^{-6}$. This can be explained as follows: the cost of evaluating $\varphi(\mathcal{J}(u)dt)f(u)$ using the Krylov projection algorithm depends on the size of the Krylov basis vector, which, in turn, depends on the step size; for a basis vector of size $m$, the cost scales as $\mathcal{O}(m^2)$. If tol $\sim 10^{-3} - 10^{-4}$, large step sizes are permitted by the traditional controller. The direct consequence of this is a large Krylov basis which turns out to be computationally prohibitive. The cost decreases with a decrease in the tolerance, i.e. a more stringent limitation on the step size may result in smaller basis vectors which may contribute to a reduction in the computational cost (here, Fig. 2). This phenomenon has also been observed by Loffeld & Tokman (2013). The proposed step-size controller does not only reduce the computational cost but, in some cases (for the Krylov method), is essential to ensure that the simulations do converge. For cases III and IV, the proposed controller improves over the performance of the traditional controller only for Leja interpolation. This is important as naively choosing an integrator (or a scheme) with a given value of tolerance may actually lead to an increase in the runtime.

Now, let us compare the performance of the Leja interpolation method with that of the Krylov subspace method. Whereas the desired systematic decrease in error with the increase in runtime is observed for Leja (all cases), the error incurred and the computational expenses are somewhat unpredictable for the Krylov subspace method (cases I and II). For this diffusion-dominated case (i.e. where the CFL condition for diffusion is more stringent than that for advection), we witness that the cost of performing a simulation using Leja interpolation is cheaper (and more accurate) than the Krylov subspace method. In case I, a comparison of the two methods based on the tolerance gives a performance enhancement of a factor of 2.5 whereas when compared with respect to (w.r.t.) the error incurred, one achieves an improvement of a factor of 1.65. In case II, this performance enhancement reaches a factor of 3. In cases III and IV (where the effects of diffusion are quite small), the Leja method is superior to the Krylov subspace method for lenient and intermediate tolerances. In terms of the error incurred, the improvement in performance of Leja over Krylov, is about a factor of 1.8 and 2.3, respectively. When compared w.r.t. tolerance, this performance enhancement reaches a factor of 6 and 8 respectively for the two cases.

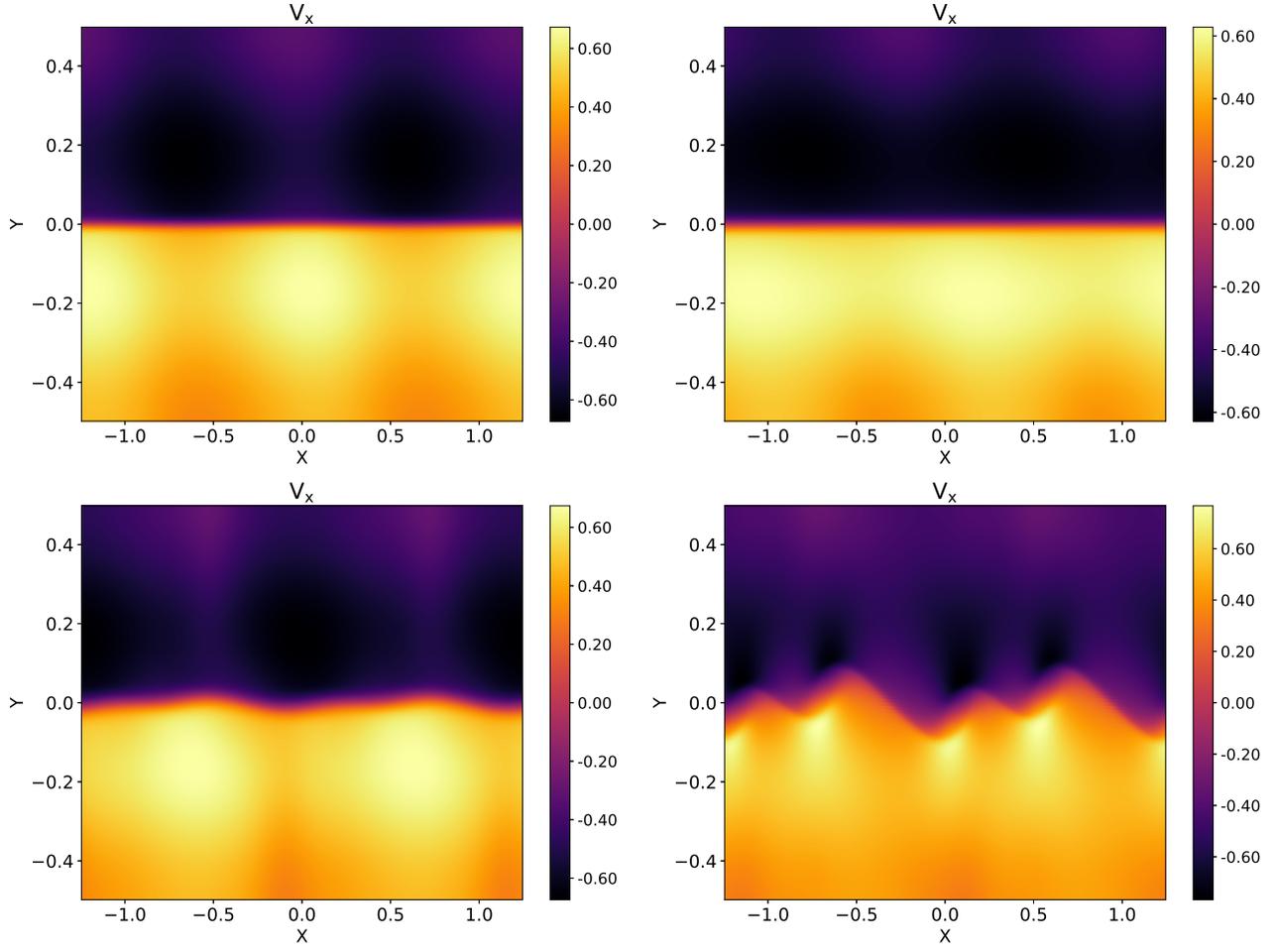

**Figure 1.** Time evolution of the velocity along the X direction at $T = 0.5$ (top-left), $T = 1.0$ (top-right), $T = 2.0$ (bottom-left), and $T = 2.75$ (bottom-right) for the KHI (case IV).

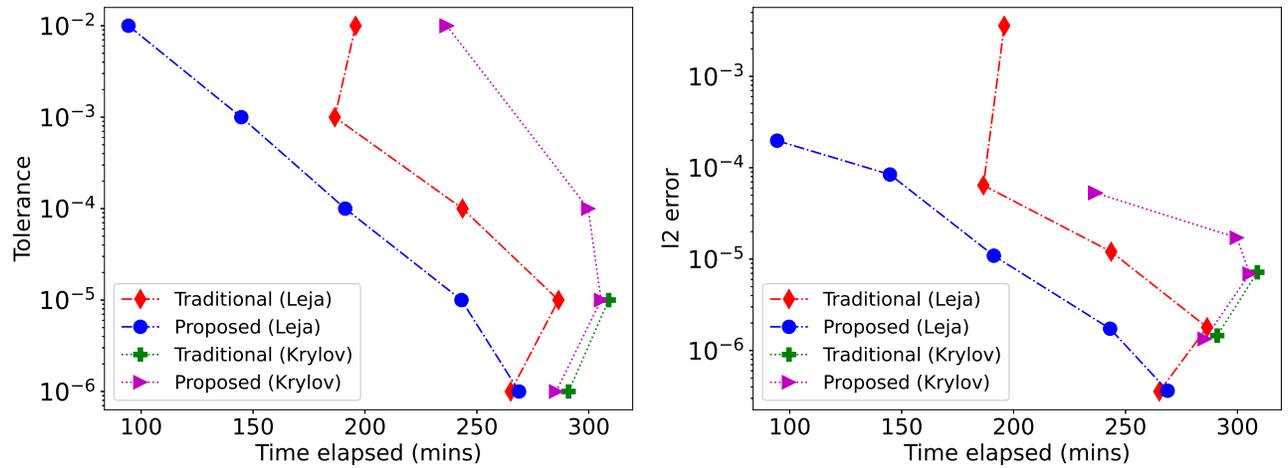

**Figure 2.** Performance of the step-size controllers for Leja interpolation and the Krylov subspace method for the EXPRB43 scheme (case I). The computational runtime is shown as a function of the tolerance (left) and the l2 norm of the global error incurred (right). The following simulations fail to converge within a reasonable amount of time: Krylov with the proposed controller for tol = $10^{-3}$ and with the traditional controller for tol $> 10^{-5}$.

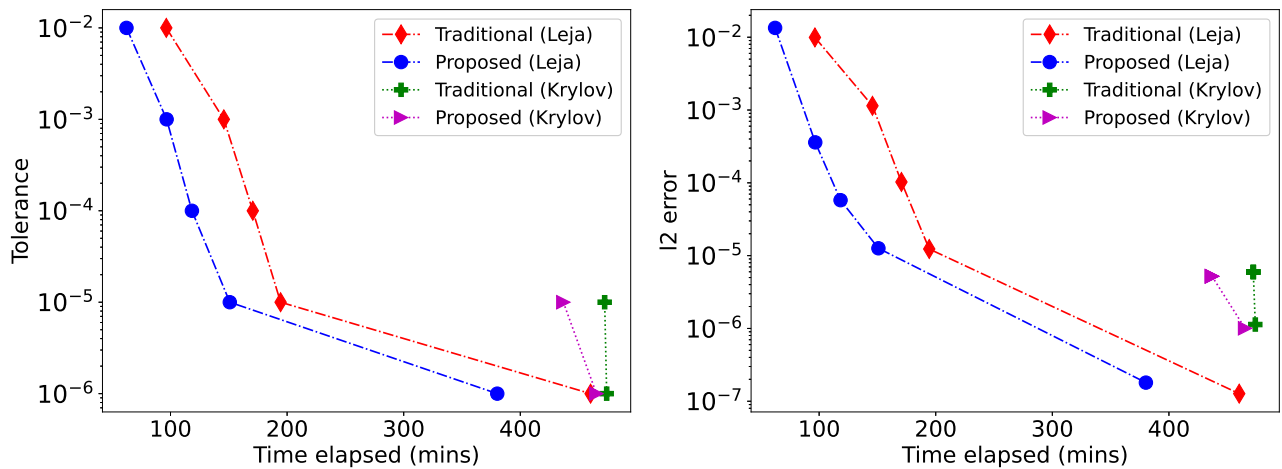

**Figure 3.** Performance of the step-size controllers for Leja interpolation and the Krylov subspace method for the EXPRB43 scheme (case II). The computational runtime is shown as a function of the tolerance (left) and the l2 norm of the global error incurred (right). The following simulations fail to converge within a reasonable amount of time: Krylov with the proposed and traditional controller for tol $> 10^{-5}$, and Leja with the proposed and traditional controller for tol $< 10^{-5}$.

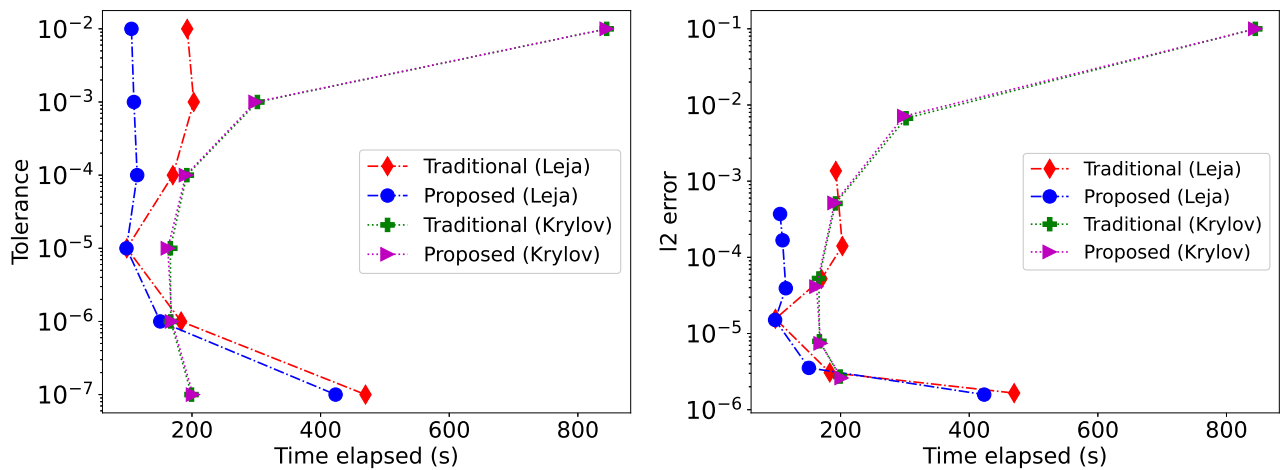

**Figure 4.** Performances of the step-size controllers for Leja interpolation and the Krylov subspace method for the EXPRB43 scheme (case III). The computational runtime is shown as a function of the tolerance (left) and the l2 norm of the global error incurred (right).

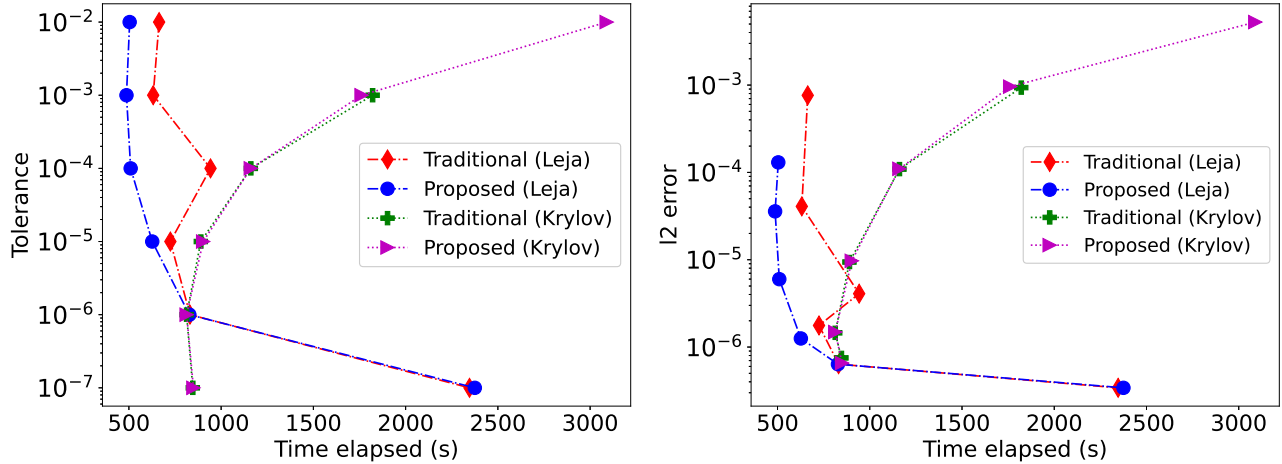

**Figure 5.** Performances of the step-size controllers for Leja interpolation and the Krylov subspace method for the EXPRB43 scheme (case IV). The computational runtime is shown as a function of the tolerance (left) and the l2 norm of the global error incurred (right). The following simulations fail to converge: Krylov with the traditional controller for tol = $10^{-3}$.

We note that the non-convergence of the Krylov-based method for lenient tolerances is due to the fact that large step sizes are permitted if the tolerance is high. This might result in errors larger than the user-defined tolerance in multiple successive time steps. Consequently, the step-size controller chooses unrealistically small step sizes to limit the error incurred in the succeeding time steps and hence the simulations fail to converge. For stringent tolerances, the chosen step size is relatively small, which is why such numerical instabilities are avoided.

#### 4.2.1. *Spectrum of the Jacobian matrix*

Interpolating a polynomial at Leja points requires some approximation of the shape of the spectrum for the matrix under consideration. In our case (matrix-free implementation in the `cppmhd` code), the most dominant eigenvalue of the relevant Jacobian is determined employing power iterations. It is not necessary to determine the exact eigenvalues at every time step. A rough estimate of the largest eigenvalue is sufficient to obtain efficient and accurate results. This, and to reduce the computational cost by quite a significant margin, is why we determine the most dominant eigenvalue for the relevant equations every 50 time steps. In Fig. 6, one can appreciate the cost reduction when the eigenvalues are determined every 50 time steps (left) as opposed to determining them at every time step (right). The total runtime of the simulations, for tol = $10^{-4}$, when the largest eigenvalue is determined every 50 time steps is around 115 s out of which 2 s (constituting ≈ 1.7% of the total runtime) are spent in determining the largest eigenvalue. This is negligible when compared to the case where the eigenvalues are determined at every time step: the total time elapsed is roughly 196 s of which 77 s (approx. 39%) are used to approximate the largest eigenvalue. This comparison has been made for five different values of the user-defined tolerance in Table 3.

Furthermore, the similarity in the shapes of the blue and red curves in Fig. 6 (right panel) clearly shows that the total computational time needed is heavily dictated by the time needed to compute the largest eigenvalue. This demonstrates that there is no advantage in determining the spectrum of the matrix at each time step. Therefore, we use this approach for all simulations that use the Leja interpolation method.

### 4.3. *Comparison of EXPRB43 with implicit and explicit schemes*

We now compare the performance of exponential integrators (EXPRB43) with that of explicit and implicit schemes. The Runge–Kutta–Fehlberg (RKF45), Dormand–Prince (DOPRI54), and Cash–Karp are some of the widely used embedded explicit schemes, each of which generates an error estimate of order four. Here, we choose the DOPRI54 scheme, and to facilitate a justified comparison with an embedded fourth-order method, we choose the RK43 (Dormand & Prince 1978) scheme as the basis for comparison with EXPRB43. In Fig. 7, we contrast the performance of EXPRB43 (Leja, blue

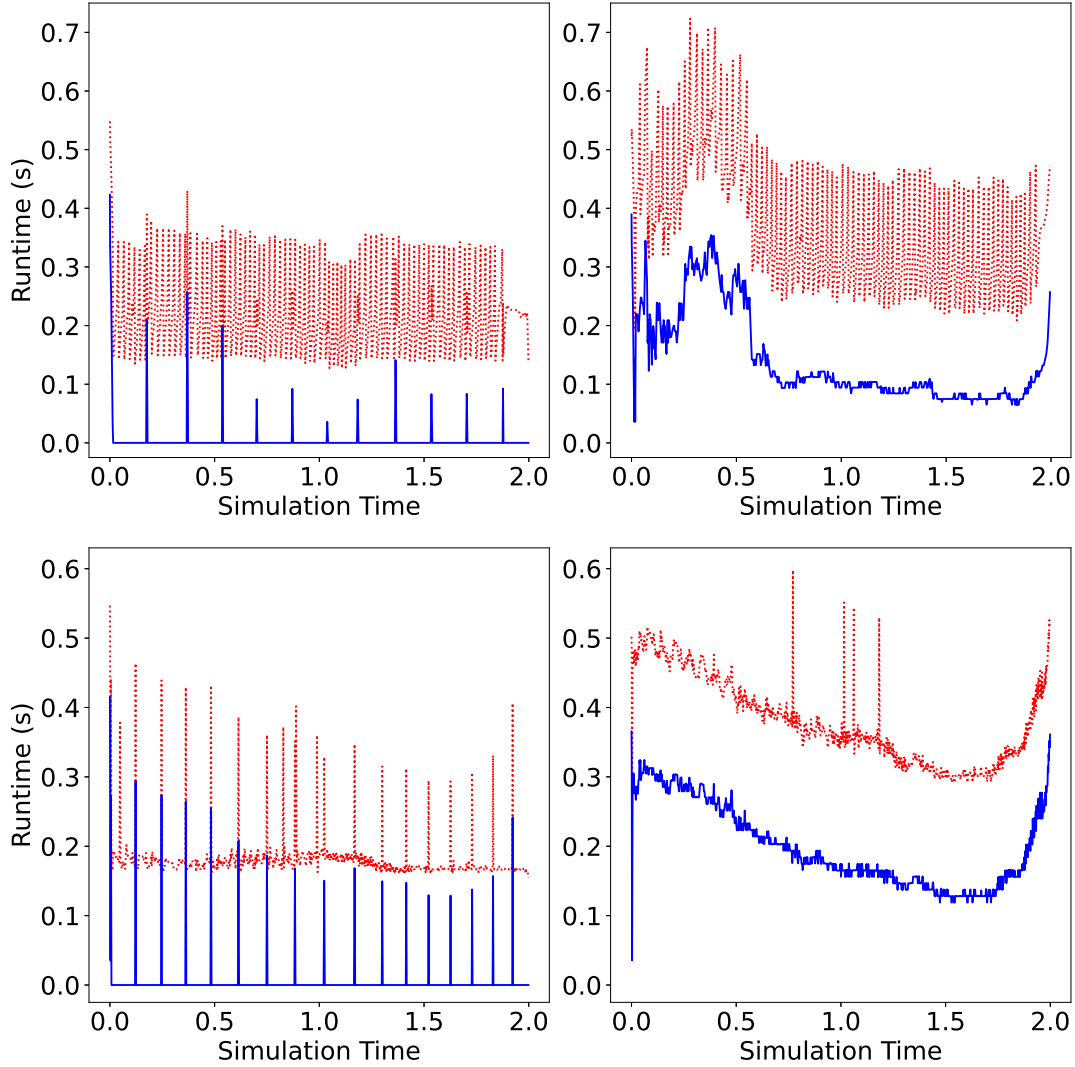

**Figure 6.** Eigenvalues determined every $50^{\text{th}}$ time step (left) vs. eigenvalues determined at every single time step (right). Red lines: total time elapsed (in seconds) during each time step, blue lines: time (in seconds) needed to compute the largest eigenvalue by means of power iterations. The simulations for this comparison (case III) are performed with the proposed controller using the EXPRB43 scheme for tolerances of $10^{-4}$ (top) and $10^{-6}$ (bottom).

**Table 3.** Cost of computing the spectrum

| Tolerance | Every 50 time steps | | Every time step | |
|---|---|---|---|---|
| | Total time (s) | Power Iters. (s) | Total time (s) | Power Iters. (s) |
| $10^{-3}$ | 110 | 2 | 165 | 63 |
| $10^{-4}$ | 115 | 2 | 196 | 77 |
| $10^{-5}$ | 98 | 3 | 211 | 108 |
| $10^{-6}$ | 151 | 4 | 320 | 167 |
| $10^{-7}$ | 426 | 19 | 1294 | 887 |

and Krylov, red) with the RK43 (purple) and DOPRI54 (green) integrator for cases I and II.

The CFL limit for diffusion (which scales as $\Delta x^2$) gives a more stringent condition compared to advection (where the CFL limit scales as $\Delta x$). As a result, the explicit schemes are severely restricted in the step sizes owing to stability constraints. This statement is well supported by the observation that even for a stringent tolerance of $10^{-6}$ ('worst-case' scenario, case II), EXPRB43 (Leja) is about 1.6 times faster than RK43 and 2.2 times faster than DOPRI54. A comparison on the basis of the l2 norm of the global error, for case I, shows that the EXPRB43 scheme, with Leja interpolation, is 3.6 and 2.4 times faster than RK43 and DOPRI54, respectively. When compared w.r.t. the user-defined tolerance, the improvement reaches a factor of 3.7 and 4.8, respectively. In case II, the performance enhancement of EXPRB43, compared to RK43 and DOPRI54, amounts to factors of 6.3 and 8.4 (based on error) and 9.8 and 13 (based on tolerance), respectively. Such a large difference in the computational cost makes an excellent case to discard the explicit schemes in favour of the exponential integrators for practical purposes.

We reiterate that the proposed step-size controller is not designed for explicit schemes (and hence the simulations with RK43 and DOPRI54 have been carried out with the traditional controller) as the computational cost per time step for these integrators is independent of the step size. This has also been illustrated with RKF45 in Deka & Einkemmer (2021).

Now, we compare the performance of EXPRB43, for cases III and IV (Fig. 8), with the variable-order backward differences formulae (BDF) implemented in the publicly available `CVODE` library. One can clearly see that for the given values of tolerance, EXPRB43 (both Leja and Krylov) is more stable in terms of the error incurred as a function of the computational cost. The performance improvement for EXPRB43 with Leja over `CVODE` reaches up to a factor of 3 and 2.3 for cases III and IV, respectively.

## 5. MAGNETIC RECONNECTION

Magnetic reconnection arises when oppositely directed magnetic field lines approach each other. This leads to a rearrangement of the field lines and conversion of a part of the energy carried by the magnetic field into kinetic energy (as well as thermal energy), resulting in the bulk of the particles being accelerated to a velocity of the order of the Alfvén velocity. The example for the reconnection problem has been drawn from Reynolds et al. (2006).

### 5.1. Simulation setup

Here, the computational domain is $X \times Y = [-12.8, 12.8] \times [-6.4, 6.4]$. The initial magnetic field distribution is given by Eq. 14. The boundary conditions are chosen to be periodic along X and reflecting along Y. Other relevant parameters are summarized in Table 4.

**Table 4.** Initial conditions for magnetic reconnection

| Parameters | Value |
| --- | --- |
| $k_x$, $k_y$ | $\pi/X$, $\pi/2Y$ |
| $\psi_0$ | 0.1 |
| $\rho$ (density) | $1.2 - \tanh^2(2y)$ |
| $P$ (pressure) | $0.5\,\rho$ |
| $(v_x, v_y, v_z)$ | $(0, 0, 0)$ |
| $\mu$ (viscosity) | $5 \cdot 10^{-2}$ |
| $\eta$ (resistivity) | $5 \cdot 10^{-3}$ |
| $\kappa$ (thermal conductivity) | $4 \cdot 10^{-2}$ |

$$\boldsymbol{B}_0 = \begin{bmatrix} \tanh(2y) - \psi_0 k_y \cos(k_x x) \sin(k_y y) \\ \psi_0 k_x \sin(k_x x) \cos(k_y y) \\ 0 \end{bmatrix} \quad (14)$$

We compare the performance of the fifth-order EPIRK5P1 with the EXPRB43 scheme for the following two cases:
**Case V:** $N_x \times N_y = 256 \times 256$; $T_f = 20$
**Case VI:** $N_x \times N_y = 128 \times 128$; $T_f = 100$
where $N_x, N_y$ correspond to the number of grid points along X and Y directions respectively, and $T_f$ is the final simulation time.

### 5.2. Results and Discussion

The work-precision diagrams for EPIRK5P1 is shown in Fig. 9 and the simulation results at the final time, for case VI, are shown in Fig. 11. The reference solution (to compute the global error) is determined with a user-specified tolerance of $10^{-11}$. The traditional as well as the proposed controller work equally well for the Krylov method for all considered cases. For the Leja method, however, the proposed controller results in an enhanced performance in some configurations whereas in others it deteriorates by a certain amount.

Comparison of the two iterative methods shows that the Krylov algorithm surpasses the Leja method for the EPIRK5P1 scheme. One of the reasons for this is that the EPIRK methods are exclusively designed to be used in conjugation with the Krylov approximation. Another reason is the caveat pointed out in Sec. 2.3. The interpolation of $\varphi(\mathcal{J}(u)\Delta t)f(u)$ as a polynomial at Leja points is highly sensitive to the magnitude of $f(u)$ and $\Delta t$. For large values of $\Delta t$ (which is naturally chosen by higher-order schemes), the polynomial fails to converge for (a maximum of) 500 Leja points. As a matter of fact, the Jacobian function diverges after just a few iterations, which forces us to choose a substantially smaller $\Delta t$. Consequently, Leja interpolation proves to be more expensive for this fifth-order integrator. To test this hypothesis, we perform further simulations of the reconnection problem with the EXPRB43 scheme, which is **not** finely tailored to be efficient with the Krylov method. The results, depicted in Fig. 10, attest to our claim that the superiority in performance of the Krylov method over Leja for the EPIRK5P1 scheme is indeed due to the "Krylov friendliness" of this scheme (see Sec. 5.3 for a detailed explanation). A comparison of the two methods for case V (EXPRB43) shows that Leja enhances the performance over Krylov by up to 87% when compared w.r.t tolerance and up to 50% when compared w.r.t. the global error incurred. In case VI, the performance enhancement for Leja reaches up to a factor of 6.4 (and 4.3) compared w.r.t tolerance (and error incurred). The results obtained with the EXPRB43 scheme clearly show the superiority of the Leja interpolation method over the Krylov method (similar to what we have seen for the KHI).

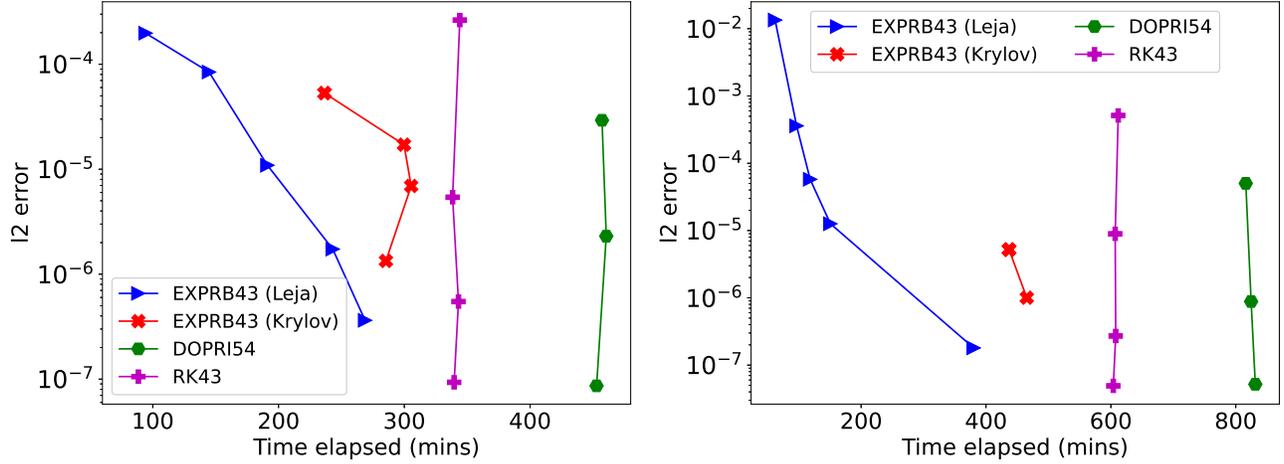

**Figure 7.** A comparison of the EXPRB43 scheme (Leja and Krylov), with the proposed controller, with the DOPRI54 scheme with the traditional controller for cases I (left) and II (right) for tol = $10^{-2}, 10^{-3}, 10^{-4}, 10^{-5}$, and $10^{-6}$. RK43 fails to converge for tol = $10^{-2}$ and DOPRI54 for $10^{-2}$ and $10^{-3}$.

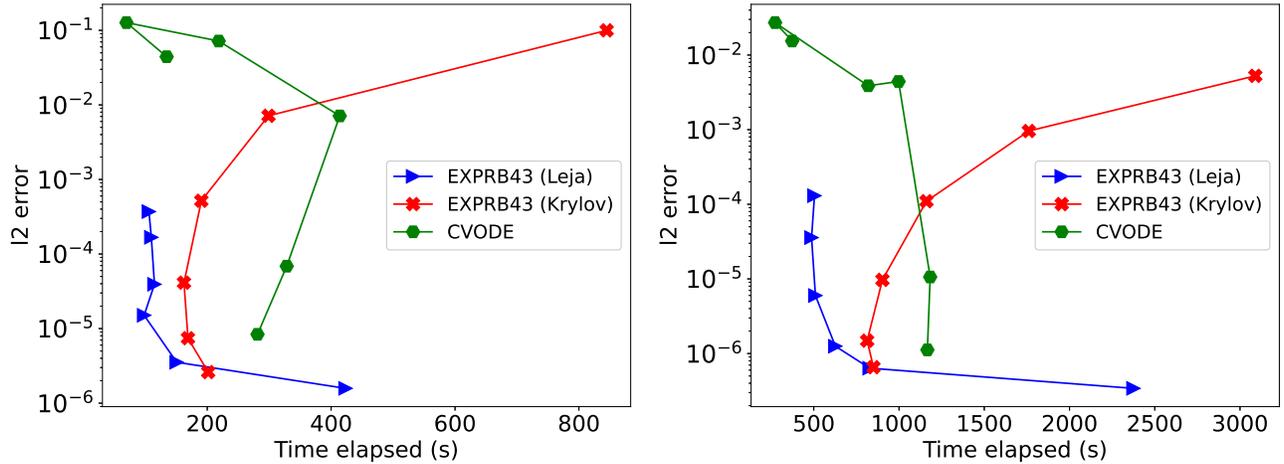

**Figure 8.** A comparison of the EXPRB43 scheme (Leja and Krylov, tol = $10^{-2}, 10^{-3}, 10^{-4}, 10^{-5}, 10^{-6}$, and $10^{-7}$), with the proposed controller, with the CVODE package (tol = $10^{-3}, 10^{-4}, 10^{-5}, 10^{-6}, 10^{-7}$, and $10^{-8}$) for cases III (left) and IV (right).

A comparison of the step-size controllers affirms that, for the Leja interpolation method, the proposed step-size controller is cheaper than the traditional controller for all considered configurations. The global error incurred by both step-size controllers is roughly the same for almost all tolerances. For stringent tolerances (tol $\leq 10^{-7}$), the step sizes become minute, leading to an increase in the total number of time steps. It is in this tolerance regime that any further decrease in the step size does not result in improved performance (in terms of computational cost), which is why both the proposed and the traditional controller are almost equally expensive for stringent tolerances. For the Krylov subspace method, the two step-size controllers have virtually the same performance.

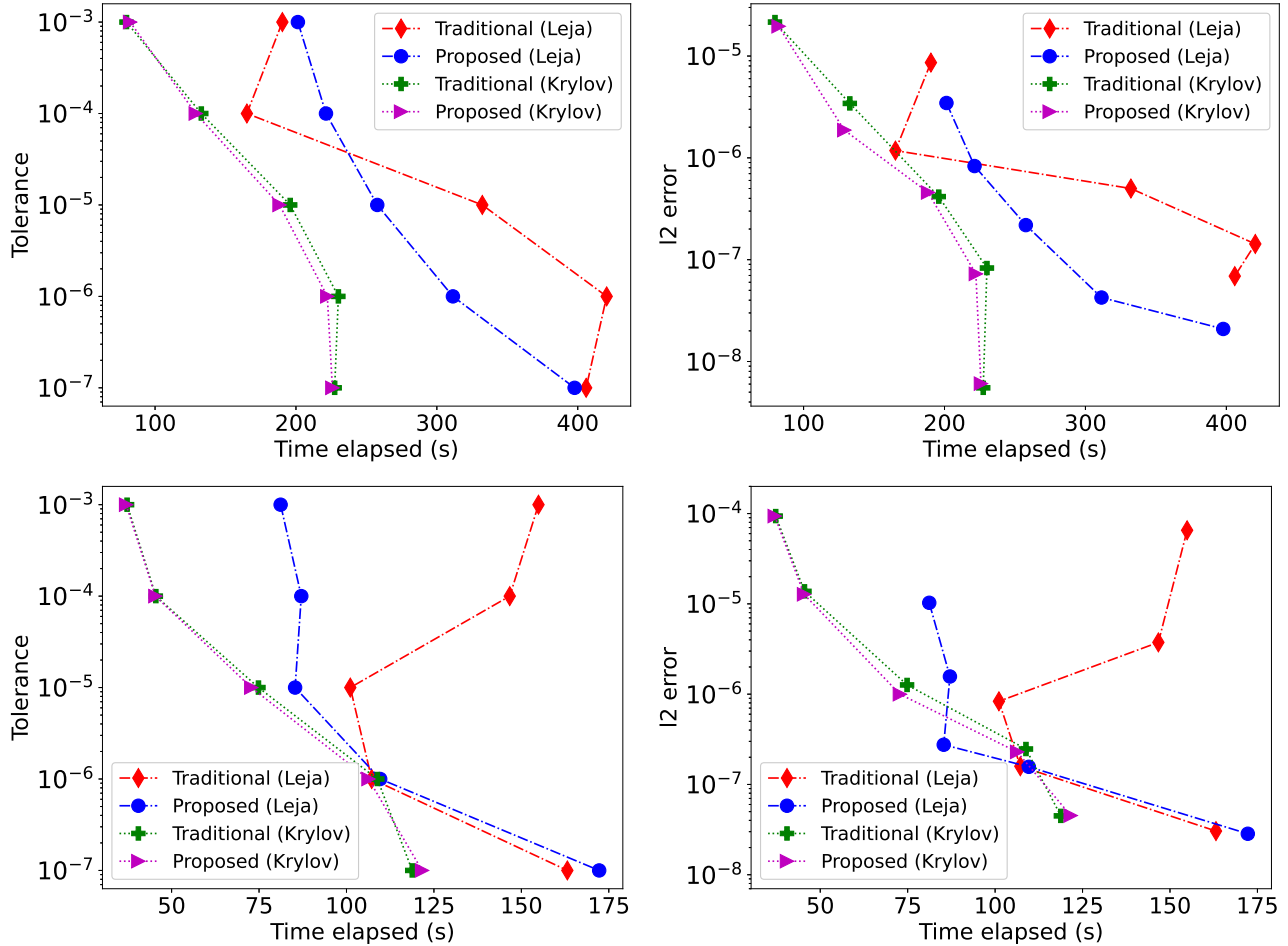

**Figure 9.** Comparison of the performances of the different step-size controllers for Leja interpolation and Krylov subspace method for the EPIRK5P1 scheme (case V - top row and case VI - bottom row). The left panel shows the computational runtime as a function of the user-defined tolerance and the right panel shows the l2 norm of the global error incurred, for the specified tolerances, as a function of the runtime.

Another noteworthy observation is that the proposed controller shows a (roughly) linear correlation between the computational cost and the tolerance (or the global error) for a large number of configurations when used with the Leja interpolation method. This is not necessarily true for the traditional controller (Figs. 9 and 10) and the Krylov method (Fig. 10).

Finally, we show that $\nabla \cdot \mathbf{B}$ is well conserved and is reasonably independent of the tolerance or the integration scheme used (for Leja interpolation, Fig. 12).

### 5.3. Leja vs. Krylov: A discussion

Generally, one would expect the Leja interpolation method to be computationally more attractive than the Krylov approximation algorithm, as no computation of inner products is needed for Leja interpolation. The results presented in the Sec. 4.2 show that the Leja method performs better than the Krylov method for EXPRB43. However, for EPIRK5P1, it is the Krylov method that dominates over Leja (Sec. 5.2). This has to do with how the EPIRK integrators are designed. EPIRK integrators are finely tailored to be compatible with the Krylov subspace method. To be more specific, the following conditions are imposed whilst developing a new EPIRK-type integrator (Tokman & Loffeld 2010):

- the ability to reuse the same basis vector to minimize the number of Krylov projections at every time step

- as higher-order $\varphi$ functions tend to converge faster, the number of Krylov vectors are to be minimized by having the higher-order $\varphi_l(z)$ functions in the subsequent stages where a new Krylov projection is to be approximated.

These so imposed conditions are redundant for optimization when used in conjunction with the Leja interpolation method. So, it is to be expected that any integrator specifically engineered to cater to these needs may not

16 Deka & Einkemmer 2021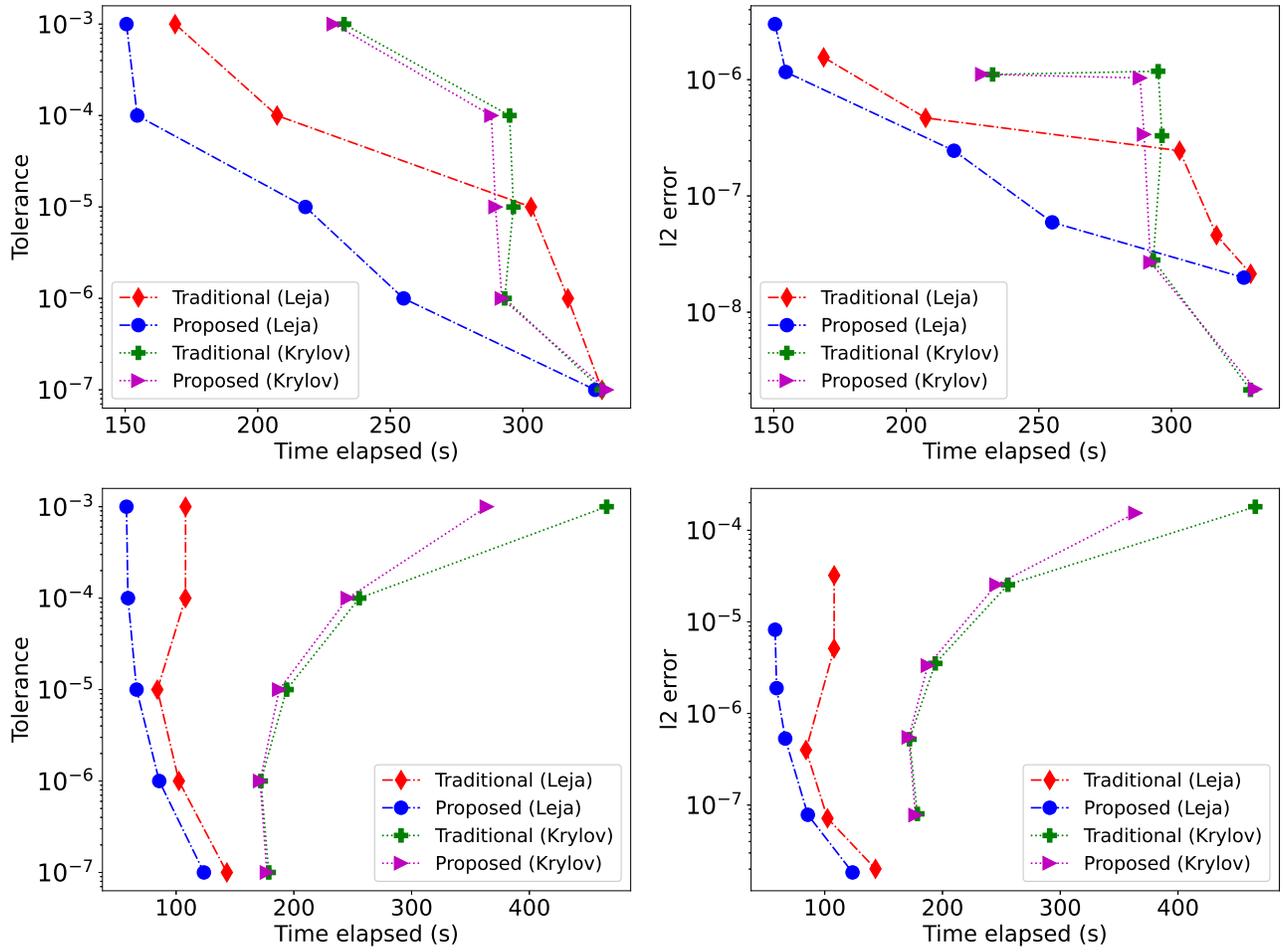

**Figure 10.** Comparison of the performance of the different step-size controllers for Leja interpolation and Krylov subspace method for the EXPRB43 scheme (case V - top row and case VI - bottom row). The left panel shows the computational runtime as a function of the user-defined tolerance and the right panel shows the l2 norm of the global error incurred, for the specified tolerances, as a function of the runtime.

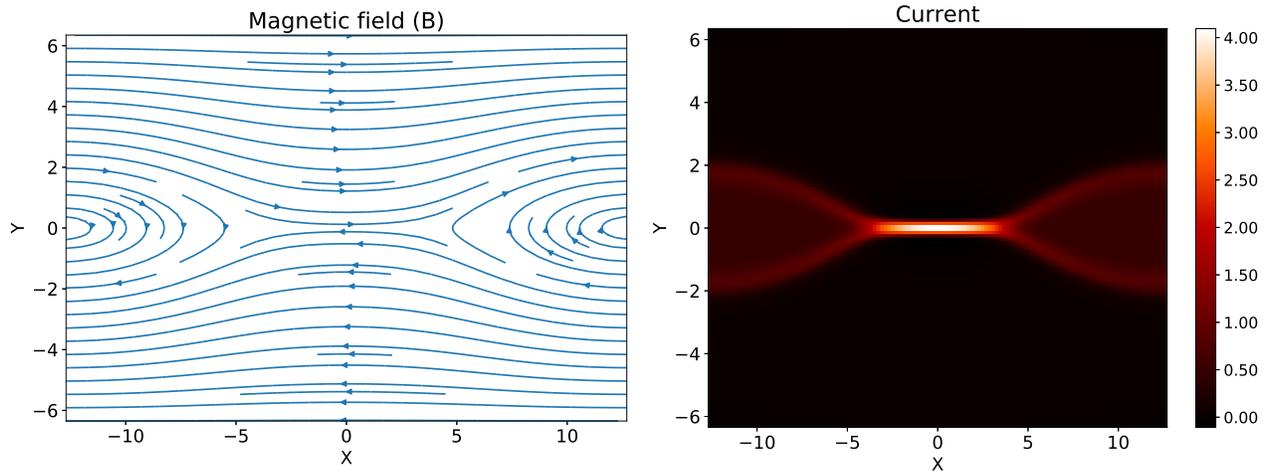

**Figure 11.** Figure shows the magnetic field (left) and the current (right) at time $T = 100$ for magnetic reconnection (case VI).

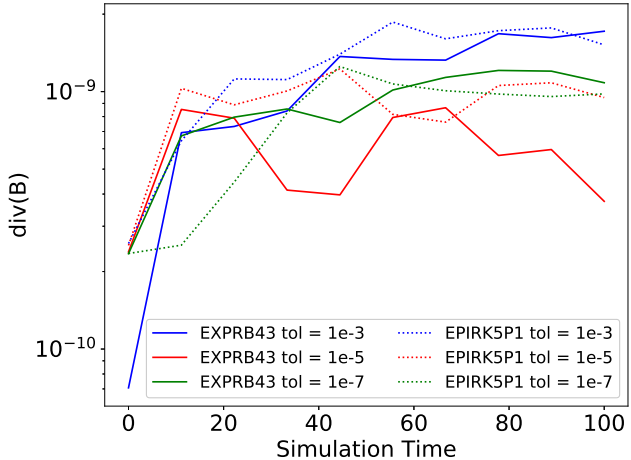

**Figure 12.** $\nabla \cdot \mathbf{B}$ is shown as a function of the simulation time for a range of tolerances (Case VI).

be well optimized for other iterative methods (e.g. the Leja method). Leja interpolation may be accelerated if multiple $\varphi_l(z)$ functions (with the same coefficients) can be grouped and applied to the same vector which may be interpolated as a polynomial. EPIRK5P1 has different coefficients (Table 1) for the internal stages, the fifth-order solution, and the error estimate. Consequently, none of the $\varphi_l(z)$ functions, so evaluated, can be reused in any of the other stages. EXPRB43 needs five matrix-function evaluations per time step, whereas EPIRK5P1 needs a total of eight matrix-function evaluations (error estimator stage included) at each time step. As such, an integrator (EPIRK5P1) finely tuned to a specific method (Krylov subspace) is expected to have predominantly improved performance.

To test if a higher-order EXPRB scheme would be competitive in the case of Leja interpolation, we consider a fifth-order Rosenbrock scheme with a fourth-order error estimate: EXPRB54s4 (Luan & Ostermann 2014). It is to be noted that this is a 4-stage integrator, and it requires ten matrix-function evaluations per time step. So, it is expected to be at least as expensive as EPIRK5P1 and certainly more than EXPRB43 (both with Leja). One might expect that this higher-order scheme should be able to take fairly large step sizes that counter the effect of a large number of matrix-function evaluations at each time step. However, as pointed out in Sec. 2.3, increasing the step size too much is not always advantageous from a performance standpoint and may cause issues with convergence. This forces one to choose a smaller step size. An interesting approach would be to rewrite a linear combination of the $\varphi_l(z)$ functions for a single stage in the form of an augmented matrix (Al-Mohy & Higham 2011), thereby reducing the total number of matrix-function evaluations. We will take this into account in our future work.

In Fig. 13, we compare the performance of this fifth-order integrator, using the Leja interpolation, with EXPRB43 and EPIRK5P1 for cases V and VI. As expected, EXPRB54s4 has worse performance than EXPRB43 for both cases. However, in case VI, EXPRB54s4 tends to marginally improve on the performance of the EPIRK5P1 scheme (with Leja), even though it needs two more matrix-function evaluations than the latter. One could speculate that a fifth-order EXPRB scheme with fewer stages and fewer matrix-function evaluations per time step could prove to be an improvement over EXPRB43 and maybe even EPIRK5P1 (with Krylov). A comprehensive analysis of the Leja and Krylov methods with several EXPRB and EPIRK integrators is needed to better understand the pros and cons of each method. This is beyond the scope of this study.

## 6. SUMMARY

We have introduced matrix-free Leja-based exponential integrators for MHD that uses an adaptive step-size controller to improve performance. We have evaluated the performance of the method for two different MHD problems. We summarize our main findings:

- Estimating the spectrum in a matrix-free implementation can be done with a negligible overhead. This has been tackled in the following way- since the most dominant eigenvalue of the Jacobian matrix essentially remains the same for several time steps, we compute the spectrum of the Jacobian only every 50 time steps. We have also shown that there is no benefit in determining the largest eigenvalue at every time step.

- Despite having to compute the spectrum of the Jacobian explicitly, Leja interpolation still outperforms the Krylov algorithm (where the size of the Krylov subspace is an implicit function of the spectrum) for the fourth-order EXPRB43 method by a significant margin.

- The proposed step-size controller has improved performance (compared to the traditional controller) for the Leja interpolation method for a wide range of user-defined tolerances. The improved performance is particularly notable in the lenient- to moderate-tolerance regime. Moreover, the zigzag curves, yielded by the traditional controller, are flattened out to a great extent. This demonstrates the practicality of using such a step-size controller with the Leja interpolation method in MHD simulations.

- The proposed controller shows a slightly enhanced performance in the Krylov subspace method for a limited number of configurations, i.e. only for the KHI. Elsewhere, its performance is comparable to that of the traditional controller (in the lenient-to the intermediate-tolerance range for magnetic reconnection).

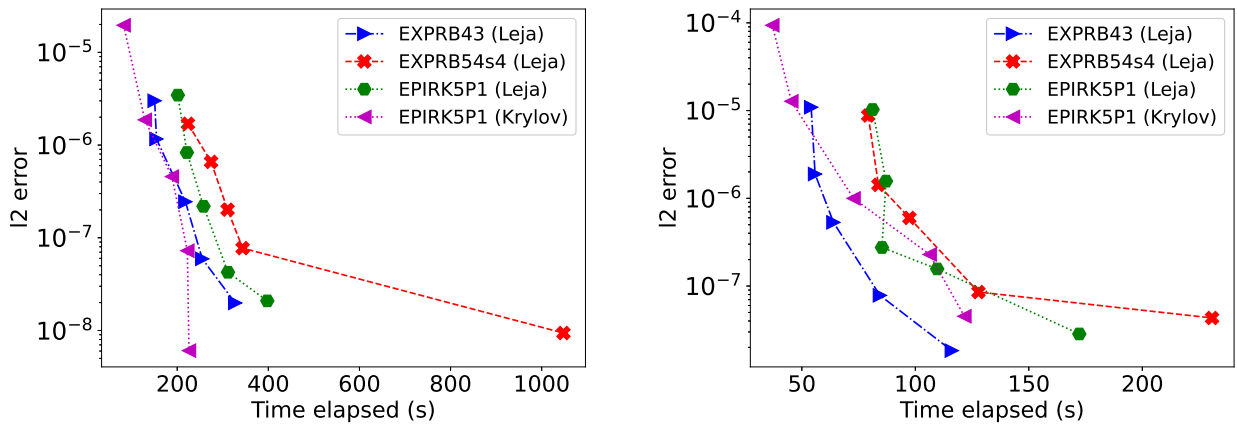

**Figure 13.** A comparison of the of the performance of the EXPRB and EPIRK schemes with Leja and Krylov for cases V and VI (left and right, respectively) with user-defined tolerances of $10^{-3}$, $10^{-4}$, $10^{-5}$, $10^{-6}$, and $10^{-7}$. All simulations have been conducted with the proposed controller.

- Multiple small step sizes often do indeed incur less computational effort than a large single step size. The proposed step-size controller can effectively exploit this fact.

- We have seen that the simulations fail to converge for large tolerances when the traditional controller with the Krylov algorithm is used (cases I and II). In these cases, the proposed controller chooses smaller step sizes, which, in turn, result in smaller Krylov basis vectors. This ensures that the simulations converge within a reasonable amount of time, which is an added benefit of the proposed controller.

- The Krylov approximation method outperforms Leja interpolation for the EPIRK5P1 scheme for almost all configurations. One of the main reasons for this is that the EPIRK schemes have been fine-tuned to be "Krylov friendly". The conditions imposed whilst designing an EPIRK scheme do not necessarily suit the needs of the Leja method. It would be interesting to construct higher-order integrators with fewer internal stages that are more favourable for Leja interpolation. This is a subject of future work.

- Embedded explicit schemes, e.g. RK43 and DO-PRI54, are not competitive with the exponential integrators. This is especially true for systems where the effects of viscosity and resistivity are substantial, i.e. problems where the explicit integrators are heavily constrained by the CFL condition.

As already mentioned, one primary advantage of the Leja interpolation method over the Krylov subspace algorithm is that it can be implemented purely in terms of matrix-vector products (neither inner products nor the explicit computation of the matrix functions in a subspace is to be considered for the Leja interpolation). In our subsequent work on the topic "Exponential Integrators for MHD," we will solve the set of resistive MHD equations on GPUs using exponential integrators in conjunction with the Leja method. This will be added as an extension to the `cppmhd` code. We also plan to study MHD problems of increasing complexity. This would potentially include the addition of more terms in the MHD equations, namely acceleration, terms involving wave-wave and wave-particle interactions, among others. Moreover, the performance of these methods in realistic scenarios will be investigated.

### ACKNOWLEDGEMENTS

This work is supported by the Austrian Science Fund (FWF)—project id: P32143-N32.

### DATA AVAILABILITY

Data related to the results presented in this work are not publicly available at the moment but can be obtained from the authors upon reasonable request.

*Software:* `cppmhd` (Einkemmer 2016), `EPIC` (Tokman 2014)

### REFERENCES

Al-Mohy, A. H., & Higham, N. J. 2011, SIAM J. Sci. Comput., 33, 488, doi: 10.1137/100788860

Almgren, A. S., Beckner, V. E., Bell, J. B., et al. 2010, ApJ, 715, 1221, doi: 10.1088/0004-637x/715/2/1221

Arnoldi, W. E. 1951, Q. Appl. Math., 9, 17, doi: 10.1090/qam/42792

Auer, N., Einkemmer, L., Kandolf, P., & Ostermann, A. 2018, Comput. Phys. Commun., 228, 115, doi: 10.1016/j.cpc.2018.02.019

Bergamaschi, L., Caliari, M., Martinez, A., & Vianello, M. 2006, in Proc. ICCS, Vol. 3994, 685–692, doi: 10.1007/11758549_93

Besse, C., Dujardin, G., & Lacroix-Violet, I. 2017, SIAM J. Numer. Anal., 55, 1387, doi: 10.1137/15M1029047

Borse, N., Acharya, S., Vaidya, B., et al. 2021, A&A, 649, A150, doi: 10.1051/0004-6361/202140440

Caliari, M., Einkemmer, L., Moriggl, A., & Ostermann, A. 2021, J. Comput. Phys., 437, 110289, doi: 10.1016/j.jcp.2021.110289

Caliari, M., Kandolf, P., Ostermann, A., & Rainer, S. 2014, BIT Numer. Math., 54, 113, doi: 10.1007/s10543-013-0446-0

Caliari, M., Vianello, M., & Bergamaschi, L. 2004, J. Comput. Appl. Math., 172, 79, doi: 10.1016/j.cam.2003.11.015

—. 2007, J. Comput. Appl. Math., 210, 56, doi: https://doi.org/10.1016/j.cam.2006.10.055

Chen, F. F. 1984, Introduction to Plasma Physics and Controlled Fusion, Vol. 1 (Springer US), doi: 10.1007/978-1-4757-5595-4

Chiuderi, C., & Velli, M. 2010, Basics of Plasma Astrophysics (Springer, Milano), doi: 10.1007/978-88-470-5280-2

Crouseilles, N., Einkemmer, L., & Prugger, M. 2018, Comput. Phys. Commun., 224, 144, doi: 10.1016/j.cpc.2017.11.003

Deka, P. J., & Einkemmer, L. 2021, arXiv:2102.02524. arxiv.org/abs/2102.02524

Dormand, J. R., & Prince, P. J. 1978, Celest. Mech., 18, 223, doi: 10.1007/BF01230162

Dreher, J., & Grauer, R. 2005, Parallel Comput., 31, 913, doi: 10.1016/j.parco.2005.04.011

Eckert, S., Baaser, H., Gross, D., & Scherf, O. 2004, Comput. Mech., 34, 377, doi: 10.1007/s00466-004-0581-1

Einkemmer, L. 2016, Comput. Phys. Commun., 206, 69, doi: 10.1016/j.cpc.2016.04.015

—. 2018, Appl. Numer. Math., 132, 182, doi: 10.1016/j.apnum.2018.06.002

Einkemmer, L., & Ostermann, A. 2015, J. Comput. Phys., 299, 716, doi: 10.1016/j.jcp.2015.07.024

Einkemmer, L., Ostermann, A., & Residori, M. 2021, Comput. Phys. Commun., 269, 108133, doi: https://doi.org/10.1016/j.cpc.2021.108133

Einkemmer, L., Tokman, M., & Loffeld, J. 2017, J. Sci. Comput., 330, 550, doi: 10.1016/j.jcp.2016.11.027

Fromang, S., Hennebelle, P., & Teyssier, R. 2006, A&A, 457, 371, doi: 10.1051/0004-6361:20065371

Goedbloed, J. P., Keppens, R., & Poedts, S. 2010, Advanced Magnetohydrodynamics: With Applications to Laboratory and Astrophysical Plasmas (Cambridge University Press), doi: 10.1017/CBO9781139195560

Gustafsson, K. 1994, ACM Trans. Math. Softw., 20, 496, doi: 10.1145/198429.198437

Gustafsson, K., Lundh, M., & Söderlind, G. 1988, BIT Numer. Math., 28, 270, doi: 10.1007/BF01934091

Hairer, E., & Wanner, G. 1996, Solving Ordinary Differential Equations II, Stiff and Differential-Algebraic Problems (Springer Berlin Heidelberg). https://link.springer.com/book/10.1007/978-3-642-05221-7

Hindmarsh, A. C., & Serban, R. 2016, User Documentation for cvode v2.9.0. https://computation.llnl.gov/sites/default/files/public/cv_guide.pdf

Hochbruck, M., & Ostermann, A. 2010, Acta Numer., 19, 209, doi: 10.1017/S0962492910000048

Hochbruck, M., Ostermann, A., & Schweitzer, J. 2009, SIAM J. Numer. Anal., 47, 786, doi: 10.1137/080717717

Huysmans, G., & Czarny, O. 2007, Nucl. Fusion, 47, 659, doi: 10.1088/0029-5515/47/7/016

Kissmann, R., Kleimann, J., Krebl, B., & Wiengarten, T. 2018, ApJS, 236, 53, doi: 10.3847/1538-4365/aabe75

Klein, C., & Roidot, K. 2011, SIAM J. Sci. Comput., 33, 3333, doi: 10.1137/100816663

Loffeld, J., & Tokman, M. 2013, J. Comput. Appl. Math., 241, 45, doi: 10.1016/j.cam.2012.09.038

Luan, V. T., & Ostermann, A. 2014, J. Comput. Appl. Math., 255, 417, doi: 10.1016/j.cam.2013.04.041

Mignone, A., Bodo, G., Massaglia, S., et al. 2007, ApJS, 170, 228, doi: 10.1086/513316

Reynolds, D. R., Samtaney, R., & Woodward, C. S. 2006, J. Comput. Phys., 219, 144, doi: 10.1016/j.jcp.2006.03.022

—. 2010, SIAM J. Sci. Comput., 32, 150, doi: 10.1137/080713331

Schilling, O. 2021, Phys. Fluids, 33, 085129, doi: 10.1063/5.0055193

Söderlind, G. 2002, Numer. Algorithms, 31, 281, doi: 10.1023/A:1021160023092

—. 2006, Appl. Numer. Math., 56, 488, doi: 10.1016/j.apnum.2005.04.026

Stone, J. M., Gardiner, T. A., Teuben, P., Hawley, J. F., & Simon, J. B. 2008, ApJS, 178, 137, doi: 10.1086/588755

Teyssier, R. 2002, A&A, 385, 337, doi: 10.1051/0004-6361:20011817

Tokman, M. 2006, J. Comput. Phys., 213, 748, doi: 10.1016/j.jcp.2005.08.032

—. 2014, EPIC (Exponential Propagation Integrators Collection). https://faculty.ucmerced.edu/mtokman/#software

Tokman, M., & Loffeld, J. 2010, Procedia Comput. Sci., 1, 229, doi: 10.1016/j.procs.2010.04.026

Tokman, M., Loffeld, J., & Tranquilli, P. 2012, SIAM J. Sci. Comput., 34, A2650, doi: 10.1137/110849961

Van Der Vorst, H. 1987, J. Comput. Appl. Math., 18, 249, doi: 10.1016/0377-0427(87)90020-3

Yang, Y., l. Wang, X., ming Li, X., mai Cao, Y., & shan Duan, W. 2021, Astrophys. Space Sci., 366, 77, doi: 10.1007/s10509-021-03984-w

Ziegler, U. 2008, Comput. Phys. Commun., 179, 227, doi: 10.1016/j.cpc.2008.02.017